\documentclass[12pt]{article}
\usepackage{amsmath,amsthm,amssymb}

\numberwithin{equation}{section}

\author{Theodore Voronov\\
\small Department of Mathematics, UC Berkeley, Berkeley, CA 92720\\
\small{\tt theodore@math.berkeley.edu}}
\title{Quantization of Forms on Cotangent Bundle}
\date{math.DG/9809130 }

\newtheorem{thm}{Theorem}[section]
\newtheorem{lm}{Lemma}[section]

\newtheorem{co}{Corollary}[section]

\theoremstyle{definition}
\newtheorem{de}{Definition}[section]
\newtheorem{ex}{Example}[section]

\newtheorem{rem}{Remark}[section]
\newtheorem*{Rem}{Remark}

\renewcommand{\a}{\alpha}

\newcommand{\e}{{\varepsilon}}

\renewcommand{\t}{{\theta}}
\newcommand{\h}{{\eta}}
\newcommand{\x}{{\xi}}
\newcommand{\z}{{\zeta}}
\renewcommand{\L}{{\Lambda}}
\renewcommand{\O}{{\Omega}}

\newcommand{\R}[1]{{\mathbb R}^{#1}}

\newcommand{\C}[1]{{\mathbb C}^{#1}}

\renewcommand{\S}[1]{\mathbb R^{0|#1}}
\newcommand{\Z}{{\mathbb Z_{2}}}
\newcommand{\ga}{{\hat \gamma^{a}}}
\newcommand{\gb}{{\hat \gamma^{b}}}
\newcommand{\Dh}{{d_{\hbar}}}

\renewcommand{\geq}{\geqslant}

\newcommand{\lder}[2]{{\partial {#1}/\partial {#2}}}
\newcommand{\der}[2]{{\frac{\partial {#1}}{\partial {#2}}}}
\newcommand{\dder}[3]{{\frac{\partial^{2} {#1}}{\partial {#2}\partial {#3}}}}
\renewcommand{\d}{{\partial}}

\newcommand{\End}{\mathop{\rm End}\nolimits}
\newcommand{\Hom}{\mathop{\rm Hom}\nolimits}
\newcommand{\Tr}{\mathop{\rm Tr}\nolimits}
\newcommand{\Str}{\mathop{\rm Str}\nolimits}
\newcommand{\tr}{\mathop{\rm tr}\nolimits}
\newcommand{\str}{\mathop{\rm str}\nolimits}
\newcommand{\Vect}{\mathop{\rm Vect}\nolimits}
\renewcommand{\div}{\mathop{\rm div}\nolimits}
\newcommand{\Ric}{\mathop{\rm Ric}\nolimits}
\newcommand{\ind}{\mathop{\rm index}\nolimits}
\newcommand{\Pf}{\mathop{\rm Pf}\nolimits}
\newcommand{\ad}{\mathop{\rm ad}\nolimits}
\newcommand{\st}{\mathbin{*_{\hbar}}}

\newcommand{\G}{\mathop{\rm \Gamma}\nolimits}

\newcommand{\const}{\frac{1}{(2\pi\hbar)^n}}
\newcommand{\sconst}{{(-i\hbar)^n}}
\newcommand{\uconst}{\frac{1}{(2\pi i)^n}}

\newcommand{\expi}[2]{{\exp_{#1}^{-1}#2}}

\newcommand{\zs}{{g_{s}(x,y)}} 

\begin{document}
\maketitle
\begin{abstract} We consider the following construction of 
quantization. For a Riemannian manifold $M$ the space of forms on 
$T^*M$ is made into a space of (full) symbols of operators acting 
on forms on $M$.  This gives rise to the composition 
of symbols, which is a deformation of the (``super'')commutative 
multiplication of forms.  The symbol calculus is exact for 
differential operators and the symbols that are polynomial in momenta.  
We calculate the symbols of natural Laplacians. (Some nice 
Weitzenb\"ock like identities appear here.)  Formulas for the traces 
corresponding to natural gradings of $\Omega (T^*M)$ are 
established.  Using these formulas, we give a simple direct 
proof of the Gauss-Bonnet-Chern Theorem. We discuss these results in the 
connection of a general question of the quantization of forms 
on a Poisson manifold.
\end{abstract}

\tableofcontents

\subsection*{Introduction}

There are two motivations for this paper. 

First, for any compact Riemannian manifold  $M$ there is a natural 
full symbol calculus for operators acting on sections of an arbitrary 
vector bundle with connection. The symbols are $\End E$-valued 
functions on $T^*M$. This symbol calculus can be viewed as a 
concrete construction of quantization for the symplectic manifold 
$T^*M$. This is more or less well-known construction (see~\cite{Bok}, 
\cite{Wid}, \cite{Ge}, \cite{Saf}, \cite{AS}, \cite{Pf}).
Suppose now that the bundle $E$ is the bundle of exterior differential 
forms. It is natural to treat forms as functions themselves 
(functions of odd variables) and then to go one step farther and 
``dequantize'' the operators from $\End \Lambda^*M$. The new symbols 
obtained in this way can be identified with  {\sl forms\/} on  
symplectic manifold $T^*M$. This is what we do in this paper.
This refined symbol calculus can be useful with the relation to the Atiyah-Singer 
Index Theorem. As an example, we give a straightforward proof of the 
Gauss-Bonnet-Chern Theorem.  One of our trace formulas 
(formula~(\ref{startrace})) can be 
utilized to prove the generalized Hirzebruch Signature Theorem.

The other motivation is as follows. Suppose we have a symplectic or a 
Poisson manifold $P$. At what conditions the Poisson structure can be 
extended from functions to differential forms~? Is it possible to quantize the 
algebra of forms~? In the present paper we provide a particular 
example of the quantization of forms, in the case when the 
manifold $P$ is a cotangent bundle $T^{*}M$. 

``Supersymmetry'' proofs of index theorems originated from the works 
of Witten  and Alvarez-Gaum\'e ~\cite{A}, who  used path integrals 
and ``physical'' considerations. First mathematically rigorous 
treatment was given for the bundle-valued Dirac operator by Getzler~\cite{Ge}, 
who combined symbol 
calculus on Riemannian manifold with Clifford algebra. This approach was 
substantially simplified 
in~\cite{AS} by an explicit use of supergeometry and deformation quantization.
However, there was some peculiarity, caused by a highly asymmetric 
way in which even and odd variables entered the 
picture (see~\cite[p.1038-1039]{AS} and Remark~\ref{getz} below).
(Symbol calculus was abandoned in~\cite{Ge2} and subsequent works in favor of the analysis 
of asymptotics of the heat kernel on the diagonal.)
In this paper we work with forms instead of spinors, and  this helps to achieve
quantization which is  more ``symmetric'' though leads to index 
formulas.

\begin{Rem}[on notation]  We use the standard language of 
supermanifolds (see~\cite{GM}, \cite{Manin}, \cite{Leites}, \cite{GIT}). 
In particular,  
and distinguish between cotangent vectors as elements of 
$T^{*}_{x}M$ and $1$-forms as elements of 
$\Pi T^{*}_{x}M$, where   $\Pi$  is the parity reversion functor. 
The same distinction is made between tangent vectors 
and ``$1$-vectors'' (multivectors of degree $1$). So on even manifold 
vectors and covectors are even, while $1$-vectors and $1$-forms are 
odd.  (As far as I know, there is no established name for the elements of 
$\Pi V$ for a vector space $V$. Borrowing some physical language, I 
can suggest to call them {\it antivectors}, and the space $\Pi V$ 
the {\it antispace}.)
\end{Rem}

\par\medskip\noindent {\bf Acknowledgments.}     The results of 
this paper were discussed at different time with M.A.~Shubin, 
A.~Weinstein, 
M.~Pflaum, J.~Rabin, L.~Sadun, S.~Rosenberg,  
and the participants of seminars at UC Berkeley, 
UCSD, University of Texas at Austin, and Boston University. 
N.Yu.~Reshetikhin attracted my attention to the problem of quantizing 
forms as differential algebra.
I am very 
much grateful to all of them.

\section{Complete symbol calculus for Riemannian manifolds}\label{QS}


Consider a Riemannian manifold  $M$ and a vector bundle  $E$ over 
$M$, with a connection. There is a natural construction that allows 
to assign pseudodifferential operators (p.d.o.'s) to functions on 
$T^*M$. More precisely, we mean p.d.o.'s acting on the sections of $E$, 
and by functions we actually mean ``$\End E$-valued functions'', i.e. 
the sections of $\pi^*(\End E)$, where $\pi: T^*M\to M$ is the 
projection. (The same can be done for 
operators $\G(E)\to \G(F)$ for two bundles, $E$, $F$.) In the following, the space of smooth 
sections of {\sl any} bundle $E$ over {\sl any} manifold $M$ is 
denoted by $C^{\infty}(M,E)$

Let us consider  the geodesic segment connecting two points  $x,y\in M$ (we 
suppose that they are sufficiently close for it to exist and be  
unique), and denote the point with the affine parameter $s$ by 
$g_{s}(x,y)$. Here $g_{0}(x,y)=x$, $g_{1}(x,y)=y$. We shall need some more 
notation. The points of $T_{x}^*M$ (cotangent vectors, or 
``momenta'', at point $x\in M$) will 
be denoted by $p$. Their coordinates will be denoted by $p_a$.  
Momenta at the point $y$ will be denoted by the letter $q$.
The parallel transport $T^{*}_{y}M\to T^{*}_{x}M$ along the geodesic, 
w.r.t. the Levi-Civita connection, will 
be denoted by $T(x,y)$. The similar transport for the bundle $E$ will 
be denoted by $\tau(x,y)$. We  shall assume that there is a 
neighborhood of the diagonal $\Delta$ in $M\times M$ where the 
geodesic segment exists and is unique. For example, let $M$ be compact 
or $M=\R{n}$. We shall use a bump function $\a(x,y)$ which is 
identically $1$ near $\Delta \subset M\times M$, nonnegative, and 
supported inside the indicated neighborhood of $\Delta$. The 
construction depends on the choice of $\a$, though not very 
heavily.

In the following we shall also need  some notation for the volume elements.
Liouville's measure on $T^*M$ is denoted by $dx\, dp$ (resp., $dy\, dq$). 
The Riemannian volume element on $M$ 
at point $x$ is $\omega(x)=\sqrt{h(x)}\, dx$.  Here 
$h(x)=\det (h_{{ab}}(x))$ is the Gram determinant in a coordinate frame.
There are Euclidean volume 
elements on $T_x M$ and $T^*_x M$. They are denoted by 
$\theta_x=\sqrt{g(x)}\,dv$ 
and $\lambda_x=(1/\sqrt{g(x)})\,dp$ respectively. Here 
$g(x)=\det (g_{ij}(x))$ is the Gram determinant in some chosen local 
frame, not necessarily coordinate. The exponential map 
$\exp_x: T_x M\to M$ allows to compare volume elements on $M$ and on 
the tangent space at a fixed point. We introduce the function 
$\mu(x,y)$ by the following equation:
\begin{equation}
  \theta_x(v)=\mu (x,y)\; \omega(y),
\end{equation}
if $y=\exp_x v$, $v\in T_x M$. This definition is valid for $y$ 
sufficiently close to $x$.

\subsection{Quantization}

\begin{de} Fix $s \in[0,1]$. Let $f\in C^{\infty}(T^*M, \End E)$. We 
associate with $f$ the following operator $\hat f$, acting on  
sections of $E$. For $u\in C^{\infty}(M, E)$
\begin{multline} \label{quant}
(\hat f u)(x):=\const\int_{T^*M}dydq \;
e^{\frac{i}{\hbar}(\expi{y}{x})\cdot q}\; \a(x,y)\,\mu(x,y)\\
 \tau(x,{\zs}) \; f({\zs},T({\zs},y)q) \;\tau({\zs},y)\;u(y).
\end{multline}
\end{de}

We can make a change of variables: $y=\exp_x v, q=T(y,x)p$. This 
yields
\begin{de}[equivalent] 
\begin{multline} \label{quantprime}
(\hat f u)(x):=\const\int_{T_x M\times T_x^*M}dvdp 
\;e^{-\frac{i}{\hbar}vp} \;\a(x,\exp_x v)\\
\tau(x,\exp_x sv) \; f(\exp_x sv,T(\exp_x sv,x)p) \;
\tau(\exp_x sv,\exp_x v)\;u(\exp_x v).
\end{multline}
\end{de}

These formulas might look complicated but the idea is very simple.

\begin{ex} Consider $M=\R{n}$. Take $\a:=1$ (no bump function is 
necessary). 
In affine coordinates on $\R{n}$ and in parallel frame for the (trivial) bundle 
$E$, the formula~(\ref{quant}) becomes
 \begin{equation} 
(\hat f u)(x):=\const\int_{\R{2n}}dydq \;e^{\frac{i}{\hbar}(y-x)\cdot 
q}\;
 f((1-s)x+sy),q) \;u(y),
\end{equation}
and the formula ~(\ref{quantprime}) becomes
\begin{equation} 
(\hat f u)(x):=\const\int_{\R{n}\times \R{n}}dvdp \; 
e^{-\frac{i}{\hbar}vp} \; f(x +sv,p) 
\;u(x + v).
\end{equation}
For $s=0$ this yields $xp$-quantization (= standard 
pseudodifferential calculus), for $s=1$ this is $px$-quantization 
(the opposite way of ordering monomials), and the case $s=1/2$ yields Weyl 
(symmetric) quantization (cp. ~\cite{BS}).
\end{ex}

\begin{ex} Take a function $f=f(x)$, not depending on $p$. Then, by 
~(\ref{quantprime}), 
\begin{multline}
	(\hat f u)(x) = \int_{T_x M}\theta_x(v) \;\delta(v)\;
\tau(x,\exp_x sv) \; f(\exp_x sv) \;
\tau(\exp_x sv,\exp_x v)\;u(\exp_x v)\\=f(x)\, u(x).
\end{multline}
(Here $\delta$ stands for the ``Euclidean'' delta-function on $T_xM$.)
\end{ex}

 \subsection{Symbols}

Let $K$ be the Schwartz kernel of $A: C^{\infty}(M,E)\to 
C^{\infty}(M,E)$ with respect to the Riemannian volume element:
\begin{equation}
	(Au)(x)=\int_M \!\omega(y) \;K(x,y)\,u(y).
\end{equation}

\begin{de} The following function on $T^*M$ is called the {\it 
symbol} of the operator $A$ and is denoted by $\sigma A$:

\begin{multline} \label{symbol}
(\sigma A)(x,p):=\int_{T_x M}\theta_x(v) \;e^{\frac{i}{\hbar}vp} 
\;\a(\exp_x(-sv),\exp_x (1-s)v)\; \rho_s(x,v)\\
\tau(x,\exp_x (-sv)) \; K(\exp_x (-sv),\exp_x (1-s)v) \;
\tau(\exp_x (1-s)v,x).
\end{multline}
Here 
\begin{equation} \rho_s(x,v):=\left( \mu(\exp_x (-sv),\exp_x 
(1-s)v)\right)^{-1}.
\end{equation}
\end{de}

Note that $\sigma A\in C^{\infty}(T^*M,\End E)$, and that $K(x,y)\in 
\Hom (E_y,E_x)$, for fixed $x,y$.

\begin{rem} For $s=0$  the symbol  of an operator $A$ can be calculated as follows:
\begin{equation}\label{specsymb}
	(\sigma A)(x,p)u_0=\left(A_y\left(e^{\frac{i}{\hbar}(\expi{x}{y})\cdot p} 
\a(x,y)\;\tau(y,x)\;u_0 \right)\right)_{|y:=x},
\end{equation}
for an arbitrary constant vector $u_0\in E_x$.  Here $A_y$ means operator 
acting on sections which argument is denoted by $y$.
\end{rem}

\begin{thm}  1. The quantization map $f\mapsto \hat f$ and the symbol 
map $A\mapsto \sigma A$ are ``almost'' mutually inverse:
\begin{align} 
	\sigma\hat f&=f(1 + O(\hbar^{\infty})), \label{sq}\\
	K_{\widehat{\sigma A}}(x,y)&=K_{A}(x,y)\cdot(\a(x,y))^2,\label{qs}
\end{align}
for any $f$ and $A$. Here $K_A$ denotes the Schwartz kernel of an 
operator $A$.

2. For polynomials in $(p_{a})$ and for differential operators the maps 
    $\ \hat{}$ and $\sigma$ are mutually inverse. In this case the 
construction is independent of the bump function $\a$.
\end{thm}

\begin{co} \label{kerdiag}If $A$ has a continuous kernel and is of 
trace class, then the traces of $A$ and $\widehat{\sigma A}$ coincide.
\end{co}
\begin{proof} By equation~(\ref{qs}), the kernels  of  $A$  and $\widehat{\sigma 
A}$ coincide  on the diagonal.
\end{proof}

\begin{rem} A complete symbol calculus for 
Riemannian manifolds was for the first time considered by 
Bokobza-Haggiag~\cite{Bok} and Widom~\cite{Wid}. Later various 
modifications were suggested and used: see~\cite{Ge}, \cite{AS}, \cite{Pf}, \cite{Saf}, 
also~\cite{Gutt}, \cite{bnw}; probably more references can be given. 
Though all approaches are 
based on the same idea (use of connection, of the exponential map and the parallel 
transport along geodesics), there are small subtleties leading to 
inequivalent constructions. In most approaches the standard 
pseudodifferential calculus (or ``$qp$''-quantization) on $\R{n}$ has 
been generalized. The beautiful idea to use a middle point on a 
geodesic segment to generalize Weyl symbols and 
other ``$s$-symbols'', compare~\cite{Shub}, is due to Yu.G.~Safarov. 
(I use the opportunity to thank M.A.~Shubin from whom I learned about Safarov's approach 
around 1989.)  The above definitions of 
quantization and symbols are based on  this idea.  For $s=0$ the construction is basically 
equivalent to the one used in~\cite{AS}, up to a slightly different 
introduction of bump function. Since we follow the philosophy 
of deformation quantization, our formulas contain Planck constant, 
which is of course very important.
\end{rem}

\subsection{Composition}

\begin{de} The operation $f\circ g:= \sigma (\hat f\hat g) $ is 
called the {\it composition } of functions $f,g\in C^{\infty}(T^*M, 
\End E)$.
\end{de}

\begin{thm} For $\hbar \to 0$,
\begin{equation} f\circ g=fg \,(1 + O(\hbar)).
\end{equation}
If $f,g$ are scalar functions, then
\begin{equation} f\circ g-g\circ f= -{i}{\hbar}\,\{f,g\}\;(1 + 
O(\hbar)),
\end{equation}
with canonical Poisson brackets on $T^*M$.
\end{thm}
\begin{proof} The first statement follows by direct computation. To 
prove the second statement, note that since in the limit $\hbar\to 0$ 
in scalar case we obtain 
commutative multiplication, the commutator 
w.r.t. the composition induces {\it some} Poisson bracket on functions. Thus it is 
sufficient to check the brackets for local coordinates $x^a,p_a$. 
That the induced brackets for them coincide with the canonical ones, follows 
from the calculation below.
\end{proof}

\subsection{Examples: quantizing linear and quadratic Hamiltonians}

In the subsequent  calculations we fix some point $x\in M  $ and work 
in normal coordinates centered at $x$. By definition, then $(\exp_x 
v)^a=v^a$, and $x^a=0$, and $\Gamma _{bc}^a(0)=0$. We also introduce a 
frame in $E$ parallel along the geodesic radii $y^a=tv^a$, which we  
call the {\it parallel frame}. It is specified by two properties:  
$A(0)=0$, $y^a A_a(y)=0$ for the connection $1$-form $A$ in $E$. This 
implies that $\d_aA_b(0)$ is antisymmetric and thus equals 
$F_{ab}(0)$, where $F$ is the curvature $2$-form. We can rewrite the 
formula~(\ref{quantprime}) in the normal coordinates and the parallel 
frame:
\begin{equation}\label{simple}
	(\hat f u)(0)=\const\int_{\R{n}\times\R{n}}dvdp 
\;e^{-\frac{i}{\hbar}vp} \;\a(x,v)\; f(sv,T(sv,0)p) \;u(v).
\end{equation}
It looks almost as for $\R{n}$.

\begin{ex} Let $X=X^a\d_a\in \Vect M$ and consider the corresponding 
fiberwise-linear Hamiltonian $X\cdot \mathrm p=X^a p_a$. By~(\ref{simple}) we have:
\begin{multline*} 
	(\widehat{X\cdot \mathrm p}\,u)(0)=\const\int_{\R{n}\times\R{n}}dvdp 
\;e^{-\frac{i}{\hbar}vp} \; X^a(sv)\,(T(sv,0)p)_a \;u(v)\\
=
\const\int_{\R{n}}dp\; p_a\int_{\R{n}}dv \;e^{-\frac{i}{\hbar}vp} \; 
(T(0,sv)X(sv))^a\;u(v)\\
=
-i\hbar \der{}{v^a}\left( (T(0,sv)X(sv))^au(v)\right)_{|v=0}\\
=-i\hbar(s\cdot \nabla_aX^a(0)+X^a(0)\,\d_a u(0)),
\end{multline*}
which we can rewrite in the invariant form:
\begin{equation}
	\widehat{X\cdot \mathrm p}=-i\hbar \,(\nabla_X+s\div X).
\end{equation}
\end{ex}

If we apply this example to $X=\d_{a}$, then we immediately obtain 
canonical commutation relations and canonical Poisson brackets.

\begin{ex}
Consider the Hamiltonian $\mathrm p^2=g^{ab}p_ap_b$, corresponding to 
the metric on $M$. In the same way,
\begin{multline*} 
	(\widehat{\mathrm p^2}u)(0)=\const\int_{\R{n}\times\R{n}}dvdp 
\;e^{-\frac{i}{\hbar}vp} \; g^{ab}(sv)\,(T(sv,0)p)_a\; 
(T(sv,0)p)_b\;u(v)\\
=\const\int_{\R{n}\times\R{n}}dvdp \;e^{-\frac{i}{\hbar}vp} \; 
g^{ab}(0)p_ap_b\;u(v)\\
=-\hbar^2 g^{ab}(0)\;\partial^2_{ab}u(0)=
-\hbar^2 g^{ab}(0)(\nabla_a\nabla_b u-\Gamma_{ba}^c\nabla_cu)(0).
\end{multline*}
We used the fact that parallel transport on $M$ is an orthogonal map. 
Thus, independent of $s$,
\begin{equation}\label{lapsymb}
	\widehat{\mathrm p^2}=-\hbar^2 \Delta,
\end{equation}
where
$\Delta$ stands for the Laplace-Beltrami operator on a Riemannian 
manifold, acting on sections of a vector bundle with a connection:
\begin{equation}
	\Delta u:=g^{ab}(\nabla_a\nabla_b\, u-\Gamma_{ba}^c\nabla_c\,u),
\end{equation}
where $\nabla_a$ in this formula denotes partial covariant derivative 
of the sections of $E$.
\end{ex}

\subsection{The trace formula}

If we take as a working definition of the trace of an operator $A$ 
the integral of its kernel over the diagonal $\Delta\subset M\times 
M$, then we easily obtain the following 
\begin{thm} \label{trace}
For all $s\in [0,1]$, the trace is calculated by the same  
formula:
\begin{align}
	\Tr \hat f&=\const\int_{T^*M} dxdp\;\tr f(x,p),\label{trq}\\
	\Tr A&=\const\int_{T^*M} dxdp\;\tr \,(\sigma A)(x,p),\label{tra}
\end{align}
for any $A$ and $f$. Here $\tr$ denotes the trace on $\End E_x$.
\end{thm}
\begin{proof} From~(\ref{quant}) we can deduce the following formula for 
the restriction of the kernel of $\hat f$ to the diagonal:
\begin{equation}
	K_{\hat f}(x,x)=\const\int_{T^*_xM}\lambda_x(p)\;f(x,p).
\end{equation}
This is true for all $s$.
Then ~(\ref{trq}) immediately follows. (Note that  
$\omega(x)\lambda_x(p)=dxdp$.) The equality~(\ref{tra}) then follows, 
by Corollary~\ref{kerdiag}.
\end{proof}

\section{Quantization of forms on $T^{*}M$}

We are going to  apply the preceding consideration to 
the particular case when $E=\Lambda(T^*M)$, the bundle of exterior 
differential  forms. (We can also consider the case of 
$\Lambda(T^*M)\otimes E$, where $E$ is  some other bundle.) In this 
case we describe endomorphisms of 
$\Lambda=\Lambda_x=\Lambda(T_x^*M)$ at each $x$ also as a result of 
quantization of some symbols. The corresponding symbol calculus ``at a 
point''  is a 
superanalog of the quantization on $\R{n}$. The calculus of the 
previous section will be combined with this pointwise calculus to 
produce our main construction.

\subsection{Quantization at a point}

Consider an arbitrary local coframe on $M$. Its 
elements are considered as $1$-forms and thus odd (see the remark in the 
introduction). So at a 
point $x$ we have $n=\dim M$ odd variables  $\xi^k=e^k(x)$. 
(Later, when we shall consider different points, we shall pick 
different letters for corresponding  values of $e^k$.) Elements of $\Lambda$ 
are functions of $\xi^k$. All endomorphisms of $\Lambda$ are 
differential operators and are generated by $\xi^k$, 
$\lder{}{\xi^k}$. They also can be expressed by integral kernels:
\begin{equation}
	(Au) (\xi)=(-1)^{n \tilde A} \int_{\S{n}} \frac{D\h}{\sqrt{g}}\; 
k(\xi,\h)\;u(\h).
\end{equation}
Here  $g=g(x)=\det 
(g_{kl}(x))$ is the Gram determinant of the frame $(e_k(x))$. Notice 
that 
\begin{equation}\label{kernel}
k_{A}(\x,\h)=(-1)^{n}(A\delta_{\h})(\x),
\end{equation}
where $\delta_{\h}(\x)=\sqrt{g}\,\delta(\x-\h)$.

\begin{lm}\label{kerprod}
\begin{equation}
k_{AB}(\x,\h) =(-1)^{n\tilde A} \int_{\S{n}} \frac{D\x'}{\sqrt{g}}\;
k_{A}(\x,\x')\;k_{B}(\x',\h)
\end{equation}
\end{lm}

In the same way as above, or, actually, in the same way as for 
$\R{n}$, we can introduce the following symbol calculus. Consider odd 
variables $\theta_k$ which transform contragrediently to $\xi^k$. 
That means that the form $\xi^k\theta_k$ is invariant under change 
of (co)frame. Geometrically, 
$\theta_k$ are the values at $x$ of the elements of the dual frame 
$e_k$ regarded as  $1$-vectors (antivectors). 
A function of $\t_k$ is a multivector at the point $x$.

In the following we fix the point $x$ and forget about it for time 
being.

\begin{de} Fix $r\in [0,1]$. 
Let $f=f(\x,\theta)$. We associate with this function the 
following operator on $\Lambda$. For any $u\in \Lambda$ set
\begin{equation} \label{squant}
	(\hat f u)(\x):=\sconst\int_{\S{2n}} 
D(\h,\t)\;e^{\frac{i}{\hbar}(\x-\h)\t}\; f((1-r)\x + r\h, \t)\;u(\h).
\end{equation}
Conversely, for an operator $A $ with the kernel $k=k(\xi,\h)$, we 
define its {\it symbol} by the formula
\begin{equation} \label{ssymb}
	(\sigma A)(\x,\t):=(-1)^{n \tilde A}\int_{\S{n}}\frac{D\e}{\sqrt{g}} 
\;e^{\frac{i}{\hbar}\e\t}\;k(\x-r\e,\x+(1-r)\e).
\end{equation}
\end{de}

\begin{rem} For $r=0$, we can express the symbol by the action of the 
operator on the exponential function (cp. ~(\ref{specsymb}):
\begin{equation}
	(\sigma A)(\x,\t)=(A_{\h}(e^{\frac{i}{\hbar}(\h-\x)\t}))_{|\h:=\x},
\end{equation}
where the subscript ${}_{\h}$ means that operator is applied to the 
functions of the variables $\h^a$.
\end{rem}

\begin{thm} The quantization map $f\mapsto \hat f$ and the symbol map 
$A\mapsto \sigma A$ are mutually inverse.
\end{thm}
\begin{proof} Since we deal with finite-dimensional spaces, it 
suffices to check either $\hat{ }\circ\!\sigma $ or 
$\sigma\!\circ\hat{ }$\/. 
Take $f=f(\x,\t)$. The integral kernel of $\hat f$, 
by~(\ref{squant}),  is
\begin{equation}\label{kerhatf}
k_{\hat f}(\x,\h)=(-1)^{n\tilde f}\sconst\int_{\S{n}} \sqrt{g}\,
D\t\;e^{\frac{i}{\hbar}(\x-\h)\t}\; f((1-r)\x + r\h, \t).
\end{equation}
So, using  formula~(\ref{ssymb}),
\begin{multline*}
\sigma(\hat f)(\x,\t)=(-1)^{n\tilde 
f}\int\frac{D\e}{\sqrt{g}}\,e^{\frac{i}{\hbar}\e\t}\;k_{\hat 
f}(\x-r\e,\x+(1-r)\e)\\=
\sconst\int\frac{D\e}{\sqrt{g}}\,\sqrt{g}\,D\t'\;e^{\frac{i}{\hbar}\e\t}
\;e^{-\frac{i}{\hbar}\e\t'}\;f((1-r)(\x-r\e)+r\x+r(1-r)\e,\t')\\=
\sconst\int D\e \,D\t' \;e^{\frac{i}{\hbar}\e(\t-\t')}\;f(\x,\t')=
\int D\t'\,\delta(\t'-\t)\,f(\x,\t')=f(\x,\t).
\end{multline*}
\end{proof}

The construction of quantization is invariant under linear 
transformations, in the following sense. Consider two spaces $V_{x}, 
V_{y}\cong \S{n}$ and an isomorphism $T: V_{y}\to V_{x}$. It induces 
pull-back of functions: $T^{*}:\Lambda_{x}\to \Lambda_{y}$. 
(Letters $x$ and $y$ are used here as just labels.)
\begin{lm} \label{spinrep}
For any  $f=f(\x_{x},\t_{x})$,
\begin{equation}
T^{*}\;\hat f\;{T^{*}}^{-1}=\widehat{T^{*}f},
\end{equation}
where $(T^{*}f)(\x_{y},\t_{y})=f(T\x_{y},T^{-1}\t_{y})$. 
\end{lm}
\begin{proof} Straightforward calculation, using the 
formula~(\ref{squant}).
\end{proof}
(This is a trivial case of the spinor representation.)
We shall use the Lemma in the next section.

The {\it composition} of symbols is defined in a usual way: $f\circ 
g:=\sigma(\hat f\hat g)$.

\begin{thm} The composition can be calculated by the following 
integral formula:
\begin{multline}\label{comp}
(f\circ 
g)(\x,\t)=\\(-i\hbar)^{2n}\!\!\!\!\!\!\!\!\int\limits_{\S{2n}\times\S{2n}}
\!\!\!\!\!\!\!\!\!
D\e_{1}\,D\tau_{1}\,D\e_{2}\,D\tau_{2}\; 
e^{\frac{i}{\hbar}(\e_{1}\tau_{2}-\e_{2}\tau_{1})}\;f(\x+r\e_{1},\t+\tau_{1})\;
g(\x+(1-r)\e_{2},\t+\tau_{2})
\end{multline}
for arbitrary $r\in [0,1]$. 
\\ For $r=0$ this formula is simplified to
\begin{equation} \label{comp0}
(f\circ g)(\x,\t)=(-i\hbar)^{n}\int_{\S{2n}}D\e\,D\tau\; e^{-\frac{i}{\hbar}
\e\tau}\;f(\x,\t+\tau)\; g(\x+\e,\t).
\end{equation}
\\ For $r=1$ it is simplified to
\begin{equation} \label{comp1}
(f\circ g)(\x,\t)=(i\hbar)^{n}\int_{\S{2n}}D\e\,D\tau\;
e^{\frac{i}{\hbar}\e\tau}\;f(\x+\e,\t)\;
g(\x,\t+\tau).
\end{equation}
 (Notice the difference in signs in the formulas~(\ref{comp0})) 
and~(\ref{comp1}).)
\end{thm}
\begin{proof}
To calculate the composition of $f$ and $g$, one has to apply twice 
the formula ~(\ref{squant}) to find $\hat f\hat g u$ for some ``test'' 
function $u$, and thus obtain the integral kernel for $\hat f\hat g$, and then 
apply~(\ref{ssymb}). After this straightforward calculation (that we 
omit), we arrive to the following formula:
\begin{multline}\label{comp00}
(f\circ g)(\x,\t)=(-i\hbar)^{2n}\int D\e\,D\t_{2}\,D\h_{1}\,D\t_{1}\;
e^{\frac{i}{\hbar}(\e\t + (\x - r\e - \h_{1})\t_{1} + (\h_{1} - 
\x - (1-r)\e)\t_{2})}\;\\
f((1-r)(\x - r\e) + r\h_{1}, \t) \; g((1-r)\h_{1} + r(\x + (1-r)\e), 
\t_{2}).
\end{multline}
Now we introduce new variables: $\e_{1}=\h_{1} - (1-r)\e  - \x$,
 $\e_{2}=\h_{1}  + r\e  - \x$. This change of variables is 
 invertible: $\h_{1}=\x + r\e_{1} + (1-r)\e_{2}$, $\e=\e_{2}-\e_{1}$, 
 and its Berezinian is unity. In these variables the formula 
 ~(\ref{comp00}) becomes
\begin{multline}
(f\circ g)(\x,\t)=(-i\hbar)^{2n}\int D\e_{1}\,D\e_{2}\,D\t_{1}\,D\t_{2}\; 
e^{\frac{i}{\hbar}(\e_{2}(\t-\t_{1}) + \e_{1}(\t_{2}-\t))}\;\\
f(\x+r\e_{1}, \t_{1})\; g(\x+(1-r)\e_{2}, \t_{2}).
\end{multline}
Finally, we substitute $\t_{1}=\t + \tau_{1}$, $\t_{2}=\t + \tau_{2}$ 
and obtain~(\ref{comp}). Formulas ~(\ref{comp0}, \ref{comp1}) follow 
from ~(\ref{comp}) by direct integration.
\end{proof}

\begin{co} For arbitrary $r\in [0,1]$ 
\begin{equation} 
(f\circ 
g)(\x,\t)=e^{-\frac{i}{\hbar}((1-r)\der{}{\t_{1}}\der{}{\x_{2}} + 
r\der{}{\x_{1}}\der{}{\t_{2}})}\;f(\x_{1},\t_{1})\;g(\x_{2},\t_{2})
\,_{\left|\parbox{2cm}{\scriptsize $ 
\begin{aligned} \x_{1}&:=\x,  & {\quad } \x_{2}&:=\x \\ 
\t_{1}&:=\t,  & { \quad} \t_{2}&:=\t \\ 
\end{aligned}$}\right.}
\end{equation}
\end{co}

\begin{co} In the limit $\hbar\to 0$ the composition becomes 
ordinary multiplication. For the commutator we have the following 
formula:
\begin{equation}
	f\circ g-(-1)^{\tilde f\tilde g}g\circ f= -i\hbar\,\{f,g\}\,(1 + O(\hbar)),
\end{equation}
where the nonvanishing Poisson brackets for the  coordinates 
$\x^a,\t_a$ are:
\begin{equation}
	\{\t_a,\x^b\}=\{\x^b,\t_a\}=\delta_a{}^b
\end{equation}
(``canonical brackets'').
\end{co}

\subsection{Formulas for local traces}
The space $\Lambda$ can be endowed with different $\Z$-gradings. In 
addition to the natural grading by parity, 
one can define the grading by the eigenvalues of an 
arbitrary  operator $S$ such that $S^2=1$. 
A form $u$ is {\it $S$-even\/} (resp., {\it $S$-odd\/}) if $Su=u$ (resp., 
$Su=-u$\/).
The {\it trivial grading\/} corresponds to $S=1$ (the identity 
operator). The operator $S$ defining the natural grading  is 
$S=(-1)^{P}$, where $P$ is the ``parity operator'', with eigenvalues 
$0, 1$. 

\begin{lm}\label{kerparity}
\begin{equation}
k_{(-1)^{P}}(\x,\h) = \sqrt{g}\delta(\x+\h)
\end{equation}
\end{lm}
\begin{proof} $k_{(-1)^{P}}(\x,\h)=(-1)^{n}((-1)^{P}\delta_{\h})(\x)=
(-1)^{n}\sqrt{g}\delta(-\x-\h)= {\sqrt{g}\delta(\x+\h)}$ 
\end{proof}

An operator $A$ that preserves the $S$-grading is called {\it 
$S$-even\/}. This is equivalent to $AS=SA$. Similarly, $A$ is {\it 
$S$-odd\/} if $AS=-SA$.

Notable example of $\Z$-grading is given by 
Hodge  ``star'' operator.  It is convenient here to present a 
``super'' construction for the star.
The Riemannian metric induces a bilinear form on the 
space $\S{n}=\Pi T^*_x M$. For $\x,\h$ denote 
$\x\cdot\h:=g_{ab}\x^a\h^b$.
We  define the operator $*$ by the following formula:
{
\renewcommand{\C}{\mathop{\rm const}\nolimits}
\begin{equation} \label{star}
	*u(\x):=\C\int_{\S{n}}\frac{D\h}{\sqrt{g}}\;e^{-it\,\x\cdot\h}\; 
u(\h).
\end{equation}
Here $t$ is an arbitrary parameter, and $\C$ stands for a normalizing factor, 
which can be 
chosen by convenience.

\begin{lm} \label{starkvadrat} In components, if $u=u_{a_1,\dots,a_k}\,\x^{a_1}\dots 
\x^{a_k}$, where we assume the coefficients skew-symmetric,   
\begin{equation}
(*u)_{a_1,\dots,a_{n-k}}=
\C\,t^{n-k}\,i^{(n-k)^2+n(n-1)}\frac{{\sqrt{g}}}{(n-k)!}\,u^{b_1,\dots,b_{k}}\;
\e_{b_1,\dots,b_k,a_1,\dots,a_{n-k}}.
\end{equation}
Here the indices are raised with the help of the metric, and 
$\e_{b_1,\dots,b_k,a_1,\dots,a_{n-k}}$ is the standard Levi-Civita 
symbol. So the operation $*$ defined by (\ref{star}) differs from  
usual only by some factor, which depends on the degree of the 
form.
For the inverse operator the following formula:
\begin{equation}
	*^{-1}=(\C)^{-2}\,(it)^{-n}(-1)^{n(n-1)/2}\;*
\end{equation}
holds for any $n$.
\end{lm}

Consider the case $n=2m$. Then  the operator $*$ is even (in the sense of natural 
grading).
Set $\C:=t^{-m}$. By Lemma~\ref{starkvadrat}, $*^2=1$. So we can consider $*$-grading. 
$*$-Even forms are what is called ``self-dual'',  and $*$-odd are 
``anti self-dual''.
}

To each choice of $\Z$-grading corresponds the {\it trace}, defined 
by the formula:
\begin{equation}
	\tr_S A:=\tr A_{00}-A_{11}
\end{equation}
for an $S$-even operator $A$, where the block decomposition 
corresponds to the chosen grading. 
Up to a factor, $\tr_S$ is specified by the property that it 
annihilates ``$S$-commutators'' (commutators of operators w.r.t. the 
chosen $\Z$-grading). 

Let us consider the traces corresponding to 
 parity, trivial grading, and $*$-grading (the last one only for 
even-dimensional case). Denote them by $\str$, $\tr_1$ and $\tr_*$ 
respectively. Obviously, $\tr_{S} A= \tr_{1} (AS)= \tr_{1} 
(SA)$ for any $S$. 

\begin{thm} \label{loctrace}
For any operator $A$, 
\begin{equation}\label{str}
	\str A= (-1)^{\tilde A n}\int\limits_{\S{n}}\frac{D\x}{\sqrt 
g}\; k_A(\x,\x)
	=(-i\hbar)^n\int\limits_{\S{2n}}\!\! D(\x,\t)\; \sigma A(\x,\t)
\end{equation}
(independent of $r$),
\begin{equation}\label{tr1}
	\tr_1 A =
	(-1)^{n\tilde A}\int\limits_{\S{n}}\frac{D\x}{\sqrt g}\; k_A(-\x,\x)
	=(-i\hbar)^n\int\limits_{\S{2n}}D(\x,\t) \;\;e^{-\frac{2i}{\hbar}\x\t}\;\sigma 
A((2r-1)\x,\t),
\end{equation}
and
\begin{equation}\label{trstar}
	\tr_* A=
	t^{-m}\!\!\!\int\limits_{\S{2m}\times\S{2m}}\frac{D\x D\h}{g}\;
	e^{-it\,\x\cdot\h}\; k_A(\x,\h)
	=t^{-m}\int\limits_{\S{2m}}\frac{D\x}{\sqrt g}\; \sigma 
	A(\x,t\hbar\,\x_{\#})
\end{equation}
(independent of $r$), where $\x_{\#a}=g_{ab} \x^{b}$.     
The last formula works for $n=2m$.
\end{thm}
\begin{proof}
Let us prove the first equality in~(\ref{str}). Consider 
$[A,B]=AB-(-1)^{\tilde A\tilde B} BA$. Applying Lemma~\ref{kerprod} 
and restricting the kernel to to the diagonal, one can deduce
$$
k_{[A,B]}(\x,\x)=\int\frac{D\x'}{\sqrt{g}}\,(-1)^{n\tilde A}\;(k_{A}(\x,\x')\,
k_{B}(\x',\x) - (-1)^{n}\,(k_{A}(\x',\x)\, k_{B}(\x,\x')),
$$
which obviously gives zero after the integration w.r.t. $\x$. Thus the 
right hand side of our formula annihilates commutators, and hence it is 
proportional to the trace $\str$. To obtain the normalization factor, it 
suffices to check some particular operator. Consider 
$A: u\mapsto u(0)$. For it, $\str A =1$ and
$k_{A}(\x,\h)=\sqrt{g}\,\h^{n}\ldots\h^{1}$, so 
$$
\int\frac{D\x}{\sqrt{g}}\;k_{A}(\x,\x)=\int 
D\x\;\x^{n}\ldots\x^{1}=1=\str A.
$$
Thus the normalization factor 
actually equals $1$. The expression for the trace via the symbol follows from the 
equality
\begin{equation}
k_{\hat f}(\x,\x)=(-1)^{n\tilde f}\sconst\int_{\S{n}} \sqrt{g}\,
D\t\; f(\x, \t)
\end{equation}
(see ~(\ref{kerhatf})). So the formula~(\ref{str})  is completely 
proved.  To prove~(\ref{tr1}), notice that $\tr_{1}A=\str ((-1)^{P}A)= 
\str (A(-1)^{P})$. By Lemmas~\ref{kerprod} and \ref{kerparity},
$$
k_{(-1)^{P}A}(\x,\h)=\int \frac{D\x'}{\sqrt{g}} 
\,\sqrt{g}\delta(\x+\x')\, k_{A}(\x',\h) = k_{A}(-\x,\h). 
$$
Applying the formula for $\str$, we obtain the first equality in 
~(\ref{tr1}). Expressing $k_{A}(-\x,\h)$ via the symbol completes the 
proof of~(\ref{tr1}). Finally, to obtain~(\ref{trstar}), we 
apply~(\ref{tr1}) to the operator $*A$. From Lemma~\ref{kerprod} and 
the formula~(\ref{star}), we obtain 
$$
k_{*A}(-\x,\x) = 
\int\frac{D\x'}{\sqrt{g}}\,t^{-m}\;e^{it\,\x\cdot\x'}\;k_{A}(\x',\x). 
$$
Thus, 
\begin{multline*}
\tr_{1}(*A)=\int\frac{D\x\,D\h}{g}\,t^{-m}\;e^{-it\x\cdot\h}\;k_{A}(\x,\h) 
=\\
\int\frac{D\x \,D\h}{g}\,t^{-m}\;e^{-it\x\cdot\h}\;
\int(-i\hbar)^{2m}\,\sqrt{g}\,D\t\;e^{\frac{i}{\hbar}(\x-\h)\t}\;
f((1-r)\x + r\h, \t).
\end{multline*}
Substitute $\x=\e-r\z$, $\h=\e + (1-r)\z$. The Berezinian of this 
change of variables is unity. Now, $(1-r)\x + r\h=\e$, $\x-\h=-\z$, 
and the argument of the first exponential becomes $-it(\e-r\z)\cdot 
(\e + (1-r)\z)=-it(\e\cdot(1-r)\z - r\z\cdot\e)= it\z\cdot\e$. So,
\begin{multline*}
\tr_{*}A=\int\frac{D\e\,D\t\,D\z}{\sqrt{g}}\,t^{-m}(-i\hbar)^{2m}\;
e^{it\z\cdot\e -\frac{i}{\hbar}\z\t}\;f(\e,\t)=\\
\int\frac{D\e\,D\t\,D\z}{\sqrt{g}}\,t^{-m}(-i\hbar)^{2m}\;
e^{-\frac{i}{\hbar}\z(\t-t\hbar\,\e_{\#})}\;f(\e,\t)=\\
\int\frac{D\e\,D\t}{\sqrt{g}}\,t^{-m}\;\delta(\t-t\hbar\,\e_{\#})\;f(\e,\t)=
\int\frac{D\e}{\sqrt{g}}\,t^{-m}\;f(\e, t\hbar\,\e_{\#}),
\end{multline*}
and the theorem is completely proved.
\end{proof}

\begin{rem} For $r=1/2$ (the Weyl quantization), 
it follows from the formula~(\ref{tr1}) that $\tr_{1} A=2^{n}\,(\sigma A) (0,0)$.
\end{rem}

\subsection{Global quantization}
Now we can combine the quantization built in Section~\ref{QS} with the 
above ``fiberwise'' quantization.  

We describe forms on $M$ as 
functions on the supermanifold $\hat M=\Pi TM$. Locally they look as 
functions $u=u(x,\x)$, where $\x^{a}$ are elements of a coframe at 
$x$. For the coordinate coframe, $\x^{a}=dx^{a}$. Operators on forms can 
be described by  ``Schwartz  kernels'' of the appearance $K(x,\x,y,\h)$:
\begin{equation}
(Au)(x,\x) = (-1)^{n\tilde A}\int_{\hat M} D(y,\h)\;K(x,\x,y,\h)\;u(y,\h),
\end{equation}
for $\x^{a}=dx^{a}$. In the
coordinate-free language, the kernel is a generalized form (de Rham's 
current) on $M\times M$, and $Au=\pi_{2*} (K\pi_{1}^{*}u)$, where 
$\pi_{1},\pi_{2}:M\times M \to M$ are projections.

After certain simplifications, we 
obtain the following formulas.  Below the letter $T$ 
stands for the parallel transport along geodesics, w.r.t. 
the Levi-Civita connection, in the bundles $T^{*}M$, $\Pi TM$ and 
$T^{*}M\Pi$.

{\sl For quantization.} Suppose we have a function $f=f(x,p,\x,\t)$.  
Geometrically, this is a section of the pull-back 
$\pi^*(\Lambda(TM\oplus T^*M))$, by $\pi: T^{*}M\to M$,   or a function 
on the supermanifold $T^*M\oplus\Pi TM\oplus T^*M\Pi)$. Denote this 
supermanifold by $N$ (this is a  bundle over $M$).
To this function we assign the following operator, acting on $\Omega 
(M)$:
\begin{multline}
	(\hat f u)(x,\x):=\uconst\int_N 
D(y,q,\h,\t)\;\a(x,y)\mu(x,y)\;e^{\frac{i}{\hbar}\left(\exp_y^{-1}x\cdot 
q + (T(y,x)\x-\h)\t\right)}\;\\
f( {\zs}, T({\zs},y)q, 
(1-r)\,T({\zs},x)\x  + r\,T({\zs},y)\h,\\  
T({\zs},y)\t )\;\;u(y,\h).
\end{multline}
Equivalent definition:
\begin{multline}\label{totquant1}
	(\hat f u)(x,\x):=\uconst\int D(v,p,\e,\t)\a(x,\exp_x 
v)\;e^{-\frac{i}{\hbar}(vp+\e\t)}\;\\
f(\exp_x sv , T(\exp_x sv,x)p,
 T(\exp_x sv,x)(\x  + r\e), T(\exp_x sv,x)\t)\\ \;u(\exp_x 
v,T(\exp_x v,x)(\x+\e)).
\end{multline}
Here integration is over the fiber $T_xM\oplus T_x^*M\oplus\Pi 
T_xM\oplus T_x^*M\Pi)$.

{\sl For symbols.} Let $A$ be an operator on differential forms, with 
the Schwartz kernel $K(x,\x,y,\h)$. 
(Here $\x^{a}$ are elements of a coframe at $x$, and 
$\h^{a}$ are elements of a coframe at $y$.)
Then its {\it symbol} $\sigma 
A=(\sigma A)\,(x,p,\x,\t)$ is defined by the formula
\begin{multline} \label{fsymbol}
	(\sigma A)(x,p,\x,\t):=\int_{T_x M\times \Pi T_x 
M}D(v,\e)\;\a(\exp_x(-sv),\exp_x (1-s)v)\; \rho_s(x,v) \\
\;e^{\frac{i}{\hbar}(vp+\e\t)} 
\;\;
K(\exp_x (-sv),T(\exp_x (-sv),x)(\x-r\e), \exp_x (1-s)v, 
\\
T(\exp_x (1-s)v,x)(\x+(1-r)\e)).
\end{multline}

These formulas can be obtained from~(\ref{quant}, \ref{symbol}) 
and~(\ref{squant}, \ref{ssymb}) by more or less direct calculation. At 
a certain step it is necessary to use Lemma~\ref{spinrep} to  change 
the points at which reference frames are taken, and to change 
variables in the integrals.

Our constructions depend on two real parameters $s,r\in [0,1]$. 
Note also that without any difficulties the definitions of 
quantization and symbols can be extended from ``scalar-valued'' forms 
to forms taking values 
in an arbitrary vector bundle with a connection. Quite helpful for 
practical calculations is the appearance of the formula~(\ref{totquant1}) 
in normal coordinates centered at $x$ (so $x^{a}=0$, 
$(\exp_{x}v)^{a}=v^{a}$):
\begin{multline}\label{totquantnorm}
	(\hat f u)(0,\x):=\uconst\int D(v,p,\e,\t)\a(x,v)\;e^{-\frac{i}{\hbar}(vp+\e\t)}\;\\
f(sv , T(sv,0) p,
 T(sv,0)(\x  + r\e), T(sv,0)\t) \;u(v,T(v,0)(\x+\e)).
\end{multline}

As it was said, the symbols in our refined symbol calculus are 
functions on the supermanifold $N=T^{*}M\oplus \Pi TM\oplus T^{*}M\Pi$. 
Roughly speaking, they are tensor products of forms and multivectors on $M$, depending 
also on a point $p$ in the cotangent space. Luckily, they can be given 
a nicer description. Using the Levi-Civita connection on $M$, the supermanifold $N$ 
can be identified with $\widehat{T^{*}M}$. That means that our 
symbols can be considered simply as forms on $T^{*}M$. For clarity, let us 
write down the change of coordinates on $N$ (assuming that we use 
coordinate frames):
\begin{align}
&\left\{ 
\begin{aligned}
x^{a}&=x^{a}(x')\\
p_{a}&=\der{x^{a'}}{x^{a}}(x(x'))\,p_{a'}\\
dx^{a}&=dx^{a'}\,\der{x^{a}}{x^{a'}}(x')\\
\t_{a}&=\der{x^{a'}}{x^{a}}(x(x'))\,\t_{a'}
\end{aligned}
\right.\\
\intertext{Compare it with  $\widehat{T^{*}M}$:}
&\left\{ 
\begin{aligned}
x^{a}&=x^{a}(x')\\
p_{a}&=\der{x^{a'}}{x^{a}}(x(x'))\,p_{a'}\\
dx^{a}&=dx^{a'}\,\der{x^{a}}{x^{a'}}(x')\\
dp_{a}&=d\left(\der{x^{a'}}{x^{a}}(x(x'))\right)\,p_{a'} + \der{x^{a'}}{x^{a}}(x(x'))\,dp_{a'}
\end{aligned}
\right.
\end{align}
The tempting idea to identify $\t_{a}$ with $dp_{a}$ fails because of 
the additional term in the transformation law for $dp_{a}$. However, 
using the connection, we can identify $\t_{a}$ with $\nabla p_{a}$, 
where
\begin{equation}
\nabla p_{a}=dp_{a} - dx^{b}\G_{ab}^{c}p_{c}.
\end{equation}
Notice that the change of variables from $(x,p,dx,dp)$ to 
$(x,p,dx,\nabla p)$ has the unit Berezinian, hence all our integrals 
over $(x,p,\x,\t)$ can be rewritten as integrals over $(x,p,dx,dp)$, 
i.e., as integrals of differential forms over the manifold $T^{*}M$.

\begin{rem} \label{getz} In the paper~\cite{Ge} Getzler proposed a 
complete symbol calculus for spinor fields, combining a symbol 
calculus on manifolds with the ``Weyl'' isomorphism $\L(V)\to C(V)$ 
(where $\L(V)$ and $C(V)$ are respectively exterior and Clifford 
algebra of a Euclidean space $V$).  The resulting symbols were {\sl 
horizontal\/} forms on $T^{*}M$, i.e., functions  $f=f(x,p,\x)$ in our notation. 
The key element of that construction  was some nonobvious
filtration introduced into the algebra of symbols (by the total degree 
in $p_{a}$ and $\x^{b}$, in our notation). 
As we showed 
in~\cite{AS}, in the quantization language this was equivalent to a 
very peculiar convention 
$\hbar_{\text{spinor}}=\hbar_{\text{usual}}^{2}$ (in~\cite{Ge} 
no Planck constant was used, so this fact was not explicit there). 
That strange convention was nevertheless essential for calculating the 
Dirac index, see ~\cite{AS}. Comparing with the complete symbol calculus 
constructed it this section, we see that now {\sl all\/} forms 
on $T^{*}M$ appear as symbols, and  odd  and  even variables are on equal 
footing (no discrimination w.r.t.  $\hbar$).
\end{rem}

\subsection{Examples}

The exterior differential and the codifferential (or divergence), 
which is defined on Riemannian manifold, are given by the 
formulas:
\begin{align}
	d&=dx^a\der{}{x^a}=\x^a\,\nabla_a  \label{d}\\
	\delta &=g^{ab}\der{}{\x^a}\,\nabla_b.\label{delta}
\end{align}
(Notice that the second expression for $d$ in ~(\ref{d}) is valid for arbitrary 
frame.) The covariant derivative of (inhomogeneous) form is
\begin{equation}\label{nabla}
\nabla_{a}=\der{}{x^{a}} - \G_{ak}{}^{l}\x^{k}\der{}{\x^{l}},
\end{equation}
where $\G_{ak}{}^{l}$ is the Christoffel symbol.

\begin{rem} In components: if $u=u_{a_1,\dots,a_k}dx^{a_1}\dots 
dx^{a_k}$, then
\begin{equation}
	(\delta u)_{a_1,\dots,a_{k-1}}=k\,u_{aa_1,\dots,a_{k-1}}{}^{;a}.
\end{equation}
\end{rem}

Direct application of the formula~(\ref{totquantnorm}) gives the 
following results for the quantization of symbols $\x\mathrm p=\x^{a} 
p_{a}$ and $\t\cdot\mathrm p=g^{ab}\t_{a} p_{b}$:

\begin{ex}
	\begin{equation}
		\widehat{\x\mathrm p}=-i\hbar\,d
	\end{equation}
\end{ex}

\begin{ex}
	\begin{equation}
		\widehat{\t\cdot\mathrm p}=-\hbar^2\,\delta
	\end{equation}
\end{ex}

\subsection{Traces}

\begin{thm}
 For any operator $A$ acting on $\Omega (M)$,
\begin{align}
\Str A &= \uconst\int_{N}D(x,p,\x,\t)\;\sigma A(x,p,\x,\t) = 
\frac{i^{n}}{(2\pi)^{n}}\int_{T^{*}M}\;\sigma A   \label{supertrace}   \\
\intertext{(integral of differential form over $T^{*}M$), independent of 
$s$ and $r$,}
\Tr_1 A &=\uconst\int_{N}D(x,p,\x,\t)\;\;e^{-\frac{2i}{\hbar}\x\t}\sigma 
A(x,p,(2r-1)\x,\t),   \label{trivtrace}
\end{align}
independent of $s$.
For $n=2m$, and for any $A$,
\begin{equation}\label{startrace}
\Tr_*A=\frac{1}{(2\pi)^{2m}}\int \frac{D(x,p,\x)}{\sqrt{g}}\;\sigma 
A(x,p,\x,\hbar^{-1}\,\x_{\#}),
\end{equation}
independent of $s$ and $r$. Here $(\x_{\#})_a=g_{ab}\x^b$. (We set in 
the last formula $t:=\hbar^{-2}$.)
\end{thm}
\begin{proof} Immediately follows from Theorem~\ref{trace} and  
Theorem~\ref{loctrace}.\end{proof}

Notice that these formulas contain no powers of $\hbar$ in front of the 
integrals, compared to~(\ref{tra}) and ~(\ref{str}, \ref{tr1}, 
\ref{trstar}). This comes naturally by cancellation of the powers of $\hbar$  
 from the ``even''
and the ``odd'' parts of our quantization. 

\section{Natural Laplacians and their symbols}

\subsection{Weitzenb\"ock formula}

For any vector bundle with a connection over Riemannian $M$ there is 
an operator
\begin{equation}
	\Delta=g^{ab}(\nabla_a\nabla_b-\Gamma_{ab}{}^c\nabla_c)
\end{equation}
(here $\nabla$ stands for the covariant derivative of the sections of 
the bundle). This Laplacian is called  Laplace-Beltrami operator 
or Bochner Laplacian. In the case when the bundle is the exterior 
bundle $\Lambda M =\Lambda (T^{*}M)$, there is another natural 
construction of a Laplacian operator, namely the Hodge Laplacian 
$\square=(d+\delta)^2=d\delta +\delta d$. The relation between these 
two operators is  given by the ``Weitzenb\"ock formula''~\cite[p.~397]{W} (see also 
~\cite{Rh}). (This name 
 is used, in general, for formulas relating any natural Laplacian 
operators in arbitrary bundle; by a ``Laplacian 
operator'' is meant an operator generalizing the usual Laplacian on 
functions in the Euclidean space, in the sense that in components  
it's  the usual Laplacian applied  componentwise, plus 
some lower-order terms.) We shall use it in the following form.

\begin{thm}\label{Weitz}
\begin{align}
	\square=(d+\delta)^2&=\Delta + \Ric_a{}^b\,dx^a\der{}{dx^b} + 
R_a{}^k{}_b{}^l\,dx^adx^b\der{}{dx^k}\der{}{dx^l}  \label{weitz1}\\
&=\Delta + \Ric_a{}^b\,dx^a\der{}{dx^b} + 
\frac{1}{2}\,R_{ab}{}^{kl}\,dx^a dx^b\der{}{dx^k}\der{}{dx^l}, \label{weitz2}
\end{align}
where $R_{ab}{}^{kl}$ is the Riemann tensor, $\Ric_a{}^b=R_{ka}{}^{kb}$ 
is the Ricci tensor, indices are raised and lowered by the Riemannian 
metric.
\end{thm}
\begin{proof} By the formulas~(\ref{d}) and (\ref{delta}),  $d+\delta=\ga\, \nabla_{a}$, 
where we define $\ga:=\x^{a}+g^{ak}\lder{}{\x^{k}}$. We need to find 
$(d+\delta)^{2}=(1/2)\,[d+\delta, d+\delta]$. Directly:
\begin{multline}
[d+\delta, d+\delta]=[\ga\, \nabla_{a},\gb\, \nabla_{b}]
\\=\ga\,[\nabla_{a},\ga]\,\nabla_{b} + [\ga,\gb]\,\nabla_{a}\nabla_{b} 
-\gb\ga\,[\nabla_{a,}\nabla_{b}]-\gb\,[\ga,\nabla_{b}]\,\nabla_{a}.
\end{multline}
Observe the following facts. First, $[\ga,\gb]=2g^{ab}$. Second, 
$[\nabla_{a,}\nabla_{b}]=-R_{abk}{}^{l}\x^{k}\lder{}{\x^{l}}$ 
(cp.(\ref{nabla})). Third, in 
the frame associated with normal coordinates centered at a given 
point $x_{0}$, the operators $\ga$ and $\nabla_{b}$ commute at 
$x_{0}$ (since the Christoffel symbols and partial derivatives of the 
metric vanish at $x_{0}$). Let's calculate in these coordinates. 
Thus, at  $x_{0}$,
$$
[d+\delta, d+\delta]=g^{ab}\,\nabla_{a}\nabla_{b}+
\gb\ga\,R_{abp}{}^{q}\x^{p}\der{}{\x^{q}}
=
2\Delta + \ga\gb\,R_{abp}{}^{q}\x^{p}\der{}{\x^{q}}.
$$
Since (by direct computation)
$$
\gb\ga=\x^{b}\x^{a}+ g^{ba} +(\x^{b}g^{ak}-\x^{a}g^{bk})\der{}{\x^{k}} + 
g^{bl}g^{ak}\der{}{\x^{l}}\der{}{\x^{k}},
$$
then, using the properties of the Riemann tensor,
\begin{multline*}
[d+\delta, d+\delta]=2\Delta - 2 \x^{a}g^{bk}\,R_{abp}{}^{q}
\der{}{\x^{k}}\x^{p}\der{}{\x^{q}}
=2\Delta -2 R_{a}{}^{k}{}_{p}{}^{q}\,\x^{a}\der{}{\x^{k}}\x^{p}\der{}{\x^{q}}
\\=2\Delta + 
2 R_{a}{}^{k}{}_{p}{}^{q}\,(\x^{a}\x^{p}\der{}{\x^{k}}\der{}{\x^{q}}
-\delta_{k}{}^{p}\x^{a}\der{}{\x^{q}})
\\=2\left(
\Delta + R_{a}{}^{k}{}_{b}{}^{l}\,\x^{a}\x^{b}\der{}{\x^{k}}\der{}{\x^{l}}
+ \Ric_{a}{}^{b}\x^{a}\der{}{\x^{b}}
\right),
\end{multline*}
and the formula~(\ref{weitz1}) is proved. To deduce~(\ref{weitz2}), one should 
raise all indices in $R_{a}{}^{k}{}_{b}{}^{l}$ (respectively lowering 
them for $\x^{a}\x^{b}$), alternate in $a,b$ (the factor $1/2$ 
appears) and apply the Ricci identity $R^{akbl}+R^{bakl}+ 
\underbrace{R^{kbal}}_{-R^{bkal}}=0$.
\end{proof}

\begin{rem} Presentation of $d+\delta$ as a ``Dirac'' operator which we 
used in the proof of Theorem~\ref{Weitz} corresponds to the action on forms 
of the Clifford algebra $C(n)$,  induced by the Riemannian metric. 
Operators $\ga$ {\sl do not\/} generate all endomorphisms of $\L$. At the 
same time, the quantization that we constructed in the previous subsection
actually realizes forms as spinors for larger Clifford algebra 
$C(n,n)\cong\End\L$, which is independent of metric. So there are two 
different Clifford module structures on the space $\L$, for two 
different Clifford algebras, which should not lead to a confusion.
\end{rem}

\begin{rem} Formulas~(\ref{weitz1}, \ref{weitz2}) can be  
also rewritten as $\square = \Delta + R$, with the operator 
$R=R_a{}^k{}_b{}^l\x^{a}\der{}{\x^{k}}\x^{b}\der{}{\x^{l}}$, 
see~\cite{Ros}. In such form the Weitzenb\"ock formula was used 
in~\cite{A}. The ``normal'' form~(\ref{weitz2}) is most 
suitable for our purposes.

\end{rem}
\subsection{The symbol of the Hodge Laplacian}

\begin{thm} For any $s\in [0,1]$,
\begin{multline} \label{symlap}
	(\sigma(\square))\,(x,p,\x,\t)=\\
=-\hbar^{-2}\,\left(\mathrm p^2 + 
\frac{1}{2}\,R_{ab}{}^{kl}\,\x^a\x^b\t_k\t_l\right) + 
i\hbar^{-1}\,(1-2r)\, \Ric_a{}^b\,\x^a\t_b + r(1-r) R
\end{multline}
\end{thm}
\begin{proof}
We apply  the
Weitzenb\"ock formula~(\ref{weitz2}).
Example~\ref{lapsymb} provides $\sigma(\Delta) = -\hbar^{-2}\, \mathrm 
p^2$, independent of $s$. All we need, is to calculate 
the symbols of the remaining terms. 
Denote $A:=\Ric_a{}^b\,\x^a\lder{}{\x^b}$,
$B:= R_{ab}{}^{kl}\,\x^a \x^b\lder{}{\x^k}\lder{}{\x^l}$. For 
simplicity we shall write below $R_a{}^b$ instead of $\Ric_a{}^b$. 
First we find the kernels:
\begin{align*}
k_{A}(\x,\h)&=R_a{}^b\x^{a}\,\der{}{\x^{b}}\,\delta(\h-\x)\,\sqrt{g},\\
k_{B}(\x,\h)&=R_{ab}{}^{kl}\,\x^{a}\x^{b}\der{}{\x^{k}}\der{}{\x^{l}}\,\delta(\h-\x)\,\sqrt{g}.
\end{align*}
Now we change variables: $\x=\z-r\e, \h=\z+(1-r)\e$, or $\e=\h-\x, 
\z=(1-r)\x+r\h$. As $\lder{}{\x}=(1-r) \lder{}{\z} - \lder{}{\e}$, it follows 
that
\begin{align}
k_{A}(\z-r\e,\z+(1-r)\e)
&=-R_a{}^b(\z^{a}-r\e^{a})\der{}{\e^{b}}\delta(\e)\,\sqrt{g}, 
\label{kA}\\
k_{B}(\z-r\e,\z+(1-r)\e)
&=R_{ab}{}^{kl}(\z^{a}-r\e^{a})(\z^{b}-r\e^{b})\dder{}{\e^{k}}{\e^{l}}
\delta(\e)\,\sqrt{g} \label{kB}.
\end{align}
Substituting~(\ref{kA}) and (\ref{kB}) into~(\ref{ssymb}) and integrating by 
parts, we obtain
\begin{multline*}
(\sigma A)(\x,\t)=\int D\e\,\delta(\e)\,\der{}{\e^{b}}
\left(
-R_a{}^b
e^{\frac{i}{\hbar}\e\t}(\x^{a}-r\e^{a})
\right)
\\=
\int D\e\,\delta(\e)\,
\left(
-R_a{}^b\,\frac{i}{\hbar}\,\t_{b}\,e^{\frac{i}{\hbar}\e\t}(\x^{a}-r\e^{a})+
R_a{}^b\,e^{\frac{i}{\hbar}\e\t}\,r\delta_{b}{}^{a}
\right) 
=\frac{i}{\hbar}\,R_a{}^b\,\x^{a}\t_{b} + rR,\\
%
\shoveleft{
(\sigma B)(\x,\t)=
\int D\e\,\delta(\e)\,\dder{}{\e^{k}}{\e^{l}}
\left(
R_{ab}{}^{kl}
e^{\frac{i}{\hbar}\e\t}(\x^{a}-r\e^{a})(\x^{b}-r\e^{b})
\right)
}
\\= \int D\e\,\delta(\e)\,\der{}{\e^{k}}
\left(
\frac{i}{\hbar}\,\t_{l}\,e^{\frac{i}{\hbar}\e\t}\,R_{ab}{}^{kl}\,
(\x^{a}-r\e^{a})(\x^{b}-r\e^{b}) + 
e^{\frac{i}{\hbar}\e\t}\,2 R_{lb}{}^{kl}(-r)(\x^{b}-r\e^{b})
\right)
\\=\int D\e\,\delta(\e)\,
\Bigl(
-{\Bigl(\frac{i}{\hbar}\Bigr)}^{2}\,\t_{l}\t_{k}\,e^{\frac{i}{\hbar}\e\t}\,
R_{ab}{}^{kl}\,(\x^{a}-r\e^{a})(\x^{b}-r\e^{b}) - 
\frac{i}{\hbar}\,\t_{l}
\,e^{\frac{i}{\hbar}\e\t}\,2 R_{kb}{}^{kl}\cdot
\\
\left.
(-r)(\x^{b}-r\e^{b})+
\frac{i}{\hbar}\,\t_{k}
\,e^{\frac{i}{\hbar}\e\t}\,2 R_{lb}{}^{kl}(-r)(\x^{b}-r\e^{b})+
e^{\frac{i}{\hbar}\e\t}\,2 R_{lk}{}^{kl}(-r)(-r)
\right)
\\=\Bigl(\frac{i}{\hbar}\Bigr)^2 R_{ab}{}^{kl}\,\x^{a}\x^{b}\t_{k}\t_{l}
+\frac{i}{\hbar}\left( -r2 R_{kb}{}^{kl}\,\x^{b}\t_{l} + r2 
 R_{lb}{}^{kl}\,\x^{b}\t_{k} \right) + r^{2}\,2R_{lk}{}^{kl}
 \\=\Bigl(\frac{i}{\hbar}\Bigr)^2 R_{ab}{}^{kl}\,\x^{a}\x^{b}\t_{k}\t_{l}
 +\frac{i}{\hbar}(-4r)\,R_{a}{}^{k}\x^{a}\t_{k} - 2r^{2} R.
\end{multline*}
(where we restored $\x$ in the argument).
Finally we obtain
\begin{multline*}
\sigma\square = \sigma\Delta + \sigma A + \frac{1}{2}\sigma B
\\=
\sigma\Delta +i\hbar^{-1}R_{a}{}^{b}\x^{a}\t_{b} + rR +
-\hbar^{2}\frac{1}{2} R_{ab}{}^{kl}\,\x^{a}\x^{b}\t_{k}\t_{l} +
i\hbar^{-1}(-2r) R_{a}{}^{b}\x^{a}\t_{b} -r^{2}R
\\=
-\hbar^{2}\left(\mathrm p^{2} +  
\frac{1}{2} R_{ab}{}^{kl}\,\x^{a}\x^{b}\t_{k}\t_{l}\right)
+i\hbar^{-1}(1-2r) R_{a}{}^{b}\x^{a}\t_{b}+r(1-r) R,
\end{multline*}
which completes the proof.
\end{proof}

\section{Gauss-Bonnet-Chern Theorem}

As it is very well known, the index of any elliptic operator 
$A: \G(M, E_{0})\to \G(M,E_{1})$ on compact manifold $M$ can be 
expressed as the supertrace of a suitable function of the corresponding ``Laplacian''
$\Delta_{A}=AA^{*}+A^{*}A : \G(M, E_{0}\oplus E_{1}) \to \G(M, E_{0}\oplus E_{1})$:
\begin{equation}
\ind A = \str \varphi(\Delta_{A}),
\end{equation}
provided the trace is absolutely convergent, and $\varphi(0)=1$. In 
particular, 
\begin{equation}
\ind A = \str e^{t\Delta_{A}},
\end{equation}
for suitable $t$. In this section we shall apply this formula 
to obtain the expression for the Euler characteristic of the 
Riemannian manifold $M$ via curvature, i.e., to deduce the 
Gauss-Bonnet-Chern Theorem. Our main tool will be the formula for the 
supertrace~(\ref{supertrace}).

Consider $D=d+\delta: \O^{ev}(M)\to \O^{od}(M)$. We have
\begin{equation}\label{euler}
\chi(M)=\ind D = \str e^{t\square}
\end{equation}
for suitable $t$. Since
\begin{equation}
\sigma e^{t\square}=e^{t\sigma(\square)}\,(1+ O(\hbar)) = 
e^{t(-\hbar^{-2}(\mathrm p^{2} +\ldots)(1 + O(\hbar))}\,(1+ 
O(\hbar))
\end{equation}
(see ~(\ref{symlap})), $t=\hbar^{2}$ seems a good choice. So
\begin{equation}
\sigma (\hbar^{2}\square)= -\mathrm p^2 - 
\frac{1}{2}\,R_{ab}{}^{kl}\,\x^a\x^b\t_k\t_l + O(\hbar),
\end{equation}
and
\begin{equation}
\sigma e^{\hbar^{2}\square}=e^{\sigma 
(\hbar^{2}\square)}\,(1+O(\hbar)) =e^{-\mathrm p^2 - 
\frac{1}{2}\,R_{ab}{}^{kl}\,\x^a\x^b\t_k\t_l} + O(\hbar).
\end{equation}
Passing to the limit $\hbar\to 0$ (which is harmless since the 
supertrace ~(\ref{euler}) is independent of $\hbar$), we obtain:
\begin{multline}
\chi(M)=\ind D=\str e^{\hbar^{2}\square}=
\uconst \int_{N} D(x,p,\x,\t)\; 
e^{-\mathrm p^2 - \frac{1}{2}\,R_{ab}{}^{kl}\,\x^a\x^b\t_k\t_l} 
\\=\begin{cases}
0, & n=2m+1\\
{\displaystyle \frac{(\sqrt{\pi})^{2m}}{(2\pi i)^{2m}}\,\int\;D(x,\x)\;
\Pf (R_{ab}{}^{kl}\,\x^a\x^b)}, & n=2m
\end{cases}
\\=
\left\{
\begin{array}{l}
   0, \qquad n=2m+1\\
   {\displaystyle (-1)^{m}\int D(x,\x)\;\Pf 
   \Bigl(%
   \frac{1}{2\pi}\underbrace{\frac{1}{2}R_{ab}{}^{kl}\,\x^a\x^b}%
                  _{\text{curvature $2$-form}}%
   \Bigr) { }}
\end{array}
\right.
=\begin{cases}
0, \qquad  n=2m+1\\
{\displaystyle\int_{M} \underbrace{\Pf \left(\frac{1}{2\pi} 
\frac{1}{2}R_{ab}{}^{kl}\,dx^a\,dx^b\right)}_{\text{Euler class of 
$TM$}} }.
\end{cases}
\end{multline}
This is the Gauss-Bonnet-Chern formula for the Euler 
characteristic.

We used the well-known formula for the 
Pfaffian as a Gaussian integral:
$$
\int_{\S{n}}D\t\;e^{-\frac{1}{2}Q^{ab}\,\t_{a}\t_{b}}=
\begin{cases}
0 & \text{for $n=2m+1$}\\ 
\Pf Q & \text{for $n=2m$}
\end{cases}
$$
where $Q^{ab}=-Q^{ba}$ (see, e.g.,~\cite{GIT}).

\begin{rem} The fact that the formula for the trace~(\ref{supertrace}) 
contains no $\hbar$  factors is crucial for our calculation. It 
comes naturally in our symbol calculus. In the calculus for spinor 
fields~\cite{AS} the similar cancellation of the powers of $\hbar$ was 
achieved only by an artifial ``asymmetric'' convention for Planck constant.
\end{rem}

\section{Discussion}
With very little effort we were able to deduce the Gauss-Bonnet-Chern 
Theorem from the trace formula~(\ref{supertrace}). Though we possess 
an analogous trace formula~(\ref{startrace}) for the $*$-grading, more subtle 
analysis of the symbol is needed in the case of the Hirzebruch 
signature operator, because 
this formula contains $\hbar^{-1}$. 

We can also discuss the results with the relation to  the  problem 
of quantization of differential algebras, in particular the algebras 
of forms. In a general setting (cp.~\cite{wein1}), 
we consider deformation of commutative 
multiplication in an algebra $A$, of the form
\begin{equation}
f\st g =fg +(-i\hbar)\,\{f,g\} + 
(-i\hbar)^{2}\,\mathop{B_{2}}\nolimits (f,g) +\ldots,
\end{equation}
starting from some Poisson bracket,
and of the differential in $A$, of the form
\begin{equation} \label{dh}
\Dh f= df + (-i\hbar)\,d_{1}f + 
(-i\hbar)^{2}\,d_{2}f+\ldots,
\end{equation}
all this being considered up to  ``gauge equivalences''
$f\mapsto f+(-i\hbar)\,a_{1}(f) +   \ldots$. (Here $d_{k}$ and $a_{k}$ 
are differential operators on $A$, and $B_{k}$ are bidifferential 
operators.)
In the first order the derivation property for $\Dh$ is equivalent to
\begin{equation}\label{der}
d\{f,g\}-\{df,g\} -(-1)^{\tilde f}\{f,dg\}
=-\Bigl(d_{1}(fg)-d_{1}f\,g -(-1)^{\tilde f}f\,d_{1}g\Bigr)
\end{equation}
The obstruction to killing $d_{1}$ or, more generally, the first 
nonvanishing term $d_{n}$ in the expansion~(\ref{dh}) , $n\geq 1$,
 by a gauge transformation is the 
cohomology class $[d_{n}]$ in the complex  of operators on $A$ with 
the differential $\ad d=[d,\phantom{d}]$. 

Consider the case of forms. We hope 
to analyze the general picture elsewhere, and here we shall make just a few 
remarks.  First, there is a question whether the deformation should 
respect the $\mathbb Z$-grading or just $\Z$-grading (parity). For the 
former option, the induced Poisson bracket should be also $\mathbb 
Z$-graded, for example the bracket of  $1$-forms should 
be $2$-form (and not a function, for example). This is highly 
unlikely. Poisson structure obtained above has the canonical brackets
\begin{equation}\label{brac}
\{p_{a},x^{b}\}=\delta_{a}{}^{b}, \quad \{\t_{a},\x^{b}\}=\delta_{a}{}^{b}
\end{equation} 
in the coordinates $(x,p,\x=dx,\t)$ and  is more complicated in 
the coordinates  $(x,p,dx,dp)$. This structure  respects parity but by no means
$\mathbb Z$-grading. The second question concerns the relation between 
the bracket and $d$. An easy computation produces the following formulas 
for the action of $d$:
\begin{align*}
dx^{a}&=\x^{a}\\
dp_{a}&=\G_{ba}^{c}\x^{b}p_{c} + \t_{a} \\
d\x^{a}&=0\\
d\t_{a}&=\G_{ba}^{c}\x^{b}\t_{c} -\frac{1}{2}R_{kla}{}^{c}\x^{k}\x^{l}p_{c}, 
\end{align*}
which resemble  Cartan's equations of structure.
Direct check shows that $d$ does not respect the  brackets~(\ref{brac}). (We 
omit the explicit formulas for $d\{f,g\}-\{df,g\} -(-1)^{\tilde 
f}\{f,dg\}$, which contain both the Christoffel symbol and the Riemann 
tensor.) 
It's not surprising, of course, because in our quantization we did not 
care about $d$. The open questions are, if $d$ can be ``adjusted'' by 
additional terms to satisfy more general condition~(\ref{der}), or if some 
other quantization for $\Omega(T^{*}M)$ can be found, which would better 
 fit into the picture incorporating $d$.


\end{document}